\documentclass[11pt]{article}
\usepackage{amssymb}
\usepackage{amsmath}
\usepackage{url}
\usepackage{graphicx}
\usepackage{wrapfig}
\usepackage[font=footnotesize]{caption}
\DeclareCaptionLabelFormat{nonumber}{#1}
\usepackage{subcaption}

\newcommand\sLP{\\[\smallskipamount]}

\newcommand\bLP{\\[\bigskipamount]}

\makeindex

\begin{document}

\title{Memories of Ian G. Macdonald}

\author
{Gert Heckman, Tom Koornwinder and Eric Opdam}

\date{}

\maketitle
\begin{abstract}
This is a slightly edited translation of a paper in Dutch which appeared in
\emph{Nieuw Archief voor Wiskunde (5)} \textbf{25} (2024), No.2, 87--90 on the
occasion of I.~G. Macdonald's death in 2023,
and aimed at a very broad mathematical audience.
First we review some of Macdonald's most important older results.
Then we focus on the period 1985--1995, when Macdonald often visited the
Netherlands
and there was much interaction between his work, notably the Macdonald
polynomials,
and the work by the authors. We end with some glimpses about Macdonald as a
person. 
\end{abstract}
\begin{wrapfigure}{R}{0.3\textwidth}
\centering
\captionsetup{width=.9\linewidth}
\includegraphics[width=0.3\textwidth]{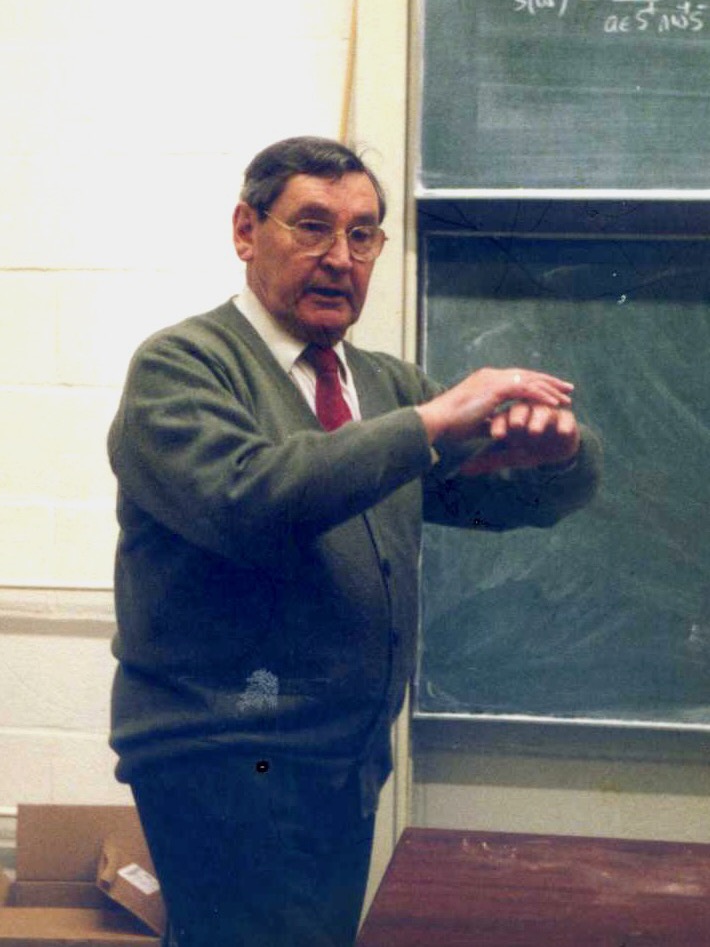}
\caption*{Macdonald giving a lecture in Amsterdam on
January 11, 2002.}
\end{wrapfigure}

Ian Grant Macdonald was an English mathematician,
known for his contributions in algebra,
geometry, Lie theory, combinatorics and special functions. He passed away on
August 8, 2023 at the age of 94.
See obituaries by Peter Cameron \cite{Cameron 2023} and
Arun Ram \cite{Ram 2024a}. Ram also wrote a longer memorial
paper~\cite{Ram 2024b}.
Much earlier A.~O.~Morris has written a nice biographical paper
\cite{Morris 2006} about him.
In the present memorial article we will in particular focus on
our contacts with Macdonald starting in the 1980s.

Macdonald studied at Trinity
College in Cambridge, where he graduated in 1952.
Instead of subsequently preparing a PhD thesis, he worked for five years as a
civil servant.
But then he could no longer resist the urge to return to academia, where he had
consecutive appointments in Manchester, Exeter, Oxford and London.
He was elected Fellow of
the Royal Society in 1979. In 2002 he received an honorary doctorate from the
University of Amsterdam. In his word of thanks he expressed his great
satisfaction to finally having received the doctor's degree.

Below we will first discuss some highlights of Macdonald's work in the 1970s.
Next we turn to his work starting in the mid 1980s on
\emph{orthogonal polynomials associated with root systems}. It was
in particular this topic where Macdonald's work had a close relationship with 
ours. There was a mutual
stimulation, by letters and by discussions during the regular visits which he
paid to our country.
What follows below is therefore not only a survey of Macdonald's work but also
a piece of history of Dutch mathematics.

\paragraph{The roots of Macdonald's work}\quad\sLP
For a better understanding of Macdonald's work one needs to have some notion
of a root system.
The theory of root systems and their symmetry groups \cite{Bourbaki 1968}
is a chapter in Euclidean geometry that can be well understood by a second year
undergraduate student in mathematics with some knowledge of linear algebra
and group theory. We give the basic definitions.

Let $V$ be a Euclidean vector space of dimension $n$ where vectors 
$\alpha,\beta$ have inner product $(\alpha,\beta)$.
For a nonzero vector $\alpha\in V$ we denote by
\[
s_{\alpha}(\lambda)=\lambda-2\,\frac{(\lambda,\alpha)}{(\alpha,\alpha)}\,\alpha 
\]
the orthogonal reflection of $\lambda$ in $V$ which uses the hyperplane
orthogonal to
$\alpha$ as a mirror. A \emph{root system} $R$ in $V$ is a finite set of
nonzero vectors in $V$ which spans $V$ and such that
\[
s_{\alpha}(\beta)\in R\quad\mbox{and}\quad
2\,\frac{(\beta,\alpha)}{(\alpha,\alpha)}\in\mathbb{Z}
\quad\mbox{for all \emph{roots} $\alpha,\beta\in R$.}
\]
Usually it is also required that for any root $\alpha$ the only multiples of
$\alpha$ in $R$ are $\pm\alpha$. Then $R$ is called a \emph{reduced} root
system. 
The subgroup $W$ of the orthogonal group $\mathrm{O}(V)$ which is generated
by the reflections $s_{\alpha}$ ($\alpha\in R$) 
is called the \emph{Weyl group} of $R$. 

The \emph{irreducible} root systems, i.e., those which cannot
be written as an orthogonal union of two smaller ones, are classified by the
symbols
\[
\mathrm{A}_n,\mathrm{B}_n,\mathrm{C}_n,\mathrm{D}_n,
\mathrm{E}_6,\mathrm{E}_7,\mathrm{E}_8,\mathrm{F}_4,\mathrm{G}_2 
\]
and accompanying Coxeter diagrams.
By the work of Claude Chevalley in the 1950s,
we know that they are also the language of the
combinatorial infrastructure of the simple algebraic 
groups $\mathbf{G}(\mathbb{F})$ over a field~$\mathbb{F}$.
The representation theory of these groups leads to
subtle identities in terms of root systems. Conversely such identities make
the representation theory computable. 

This was first understood around 1925 by Hermann Weyl, who gave  
an explicit formula, in terms of roots, for the irreducible characters 
of a compact simple Lie group, considered as a
compact real form of $\mathbf{G}(\mathbb{C})$. Weyl proved his formula by
transcendental methods.
Freudenthal gave in 1954 an algebraic proof of the character formula by means
of the quadratic  \emph{Casimir operator}. 

Something similar occurred in 1985 in Harish-Chandra's description
\cite{Harish-Chandra 1958} of the
elementary spherical functions on a Riemannian symmetric space
$\mathbf{G}(\mathbb{R})/K$, where $K$ is a maximally compact subgroup
of the non-compact simple real Lie group $\mathbf{G}(\mathbb{R})$.
These functions are eigenfunctions of a commuting system of linear
partial differential operators. Here again the explicit expression, in terms of
roots, of the quadratic Casimir operator plays a dominant role.

\paragraph{Important results of Macdonald}\quad\sLP
In the area just discussed Macdonald has obtained essential results. 
We present three of them in chronological order.

The first result is the $p$-adic analogue of the elementary spherical functions
on $\mathbf{G}(\mathbb{Q}_p)/\mathbf{G}(\mathbb{Z}_p)$,
announced by Macdonald \cite{Macdonald 1968} in 1968.
The details appeared in a book \cite{Macdonald 1971} which is based
on lectures at the Ramanujan Institute in Madras in 1970 and which is written
in a very accessible way. In the old days,
one could get a free copy by just sending a request to the author.
In a sense, the $p$-adic case is simpler than Harish-Chandra's real case, 
because the elementary spherical functions for the $p$-adics are given
by an explicit formula (for $\mathbf{G}=\mathbf{GL}_n$ in terms of
Hall--Littlewood polynomials \cite[\S V.3]{Macdonald 1995}).
Macdonald's explicit formula was one of the ingredients used by
Robert Langlands for the formulation of the conjectures named after him
\cite[p.39]{Langlands 1970}.

For a good understanding of the second result \cite{Macdonald 1972}
we first say somewhat more about
\emph{Weyl's character formula}. This gives the value of the character on a
maximal torus as a certain quotient. The denominator is a Fourier polynomial,
skew-symmetric under the Weyl group $W$ and containing a parameter $\lambda$
(\'{E}lie Cartan's \emph{highest weight}) which characterizes the irreducible 
representation.
The denominator is the famous \emph{Weyl denominator}
$\Delta=\prod_{\alpha>0} (e^{\alpha/2}-e^{-\alpha/2})$,
where the product runs over the set of positive roots $\alpha>0$,
i.e., the roots
lying on one side of a fixed generic linear subspace of $V$ of codimension~1.
When we apply the character formula to the trivial representation, the quotient
expresses the trivial character $1$. 
Hence the numerator specialised at $\lambda=0$ equals the Weyl denominator
$\Delta$, which is the celebrated \emph{Weyl denominator formula}
\[
\prod_{\alpha>0} (e^{\alpha/2}-e^{-\alpha/2})=
\sum_{w\in W} \epsilon(w) e^{w(\rho)},
\]
where $\epsilon(w)=\det(w)$ and $\rho$ is half the sum of the positive roots.
In this formula of the form $\prod=\sum$
the product has 
terms with opposite signs which cancel each other, so that fewer terms remain
in the sum than first thought.

Macdonald observed a similar phenomenon in \emph{Jacobi's
triple product identity}
\[
\sum_{m\in\mathbb{Z}}(-1)^mq^{m(m-1)/2}z^m=
\prod_{n=1}^{\infty}(1-q^n)(1-q^nz^{-1})(1-q^{n-1}z).
\]
This is a classical formula in terms of theta functions in which $q,z$
are complex variables (with $z\ne0$ and $|q|<1$).
In this formula Macdonald recognized the structure of the affine root system
of type~$\mathrm{A}_1$. Next he found
analogous identities for all other affine root systems~\cite{Macdonald 1972}.
By specialization he also obtained simple summation formulas 
for $\eta^d$, where \emph{Dedekind's eta function} is given by
\[
\eta=q^{1/24}\prod_{n=1}^{\infty} (1-q^n)
\]
and $d$ is the dimension of a compact
simple Lie group. For some values of $d$ this identity was earlier found by
the theoretical physicist Freeman Dyson.
In an interesting paper \cite{Dyson 1972} he regretted his missed opportunity
because the physicist Dyson and the mathematician Dyson had not talked with
each other. Alternative proofs of the
Macdonald identities were later given by
Victor Kac~\cite{Kac 1974} and Eduard Looijenga~\cite{Looijenga 1976}.

Macdonald felt a special sympathy for the general
linear group $\mathbf{GL}_n$ and the symmetric group $S_n$, 
the ``All-embracing Majesties'' of the linear algebraic 
and the finite groups, respectively (as Weyl \cite[end of Ch.~IV]{Weyl 1939}
phrased it for $\mathbf{GL}_n$ in ironic adoration).
Here the representation theory has a strongly combinatorial
flavour, with methods which are, in spite of their elementary nature, not
necessarily simple. In 1979 he published the influential book
\emph{Symmetric Functions and Hall Polynomials} on this matter.
A much expanded second edition  \cite{Macdonald 1995}
appeared in 1995. For this work he received
in 2009 the Leroy P. Steele Prize for Mathematical Exposition.

\paragraph{Macdonald's contact with the Low Countries}\quad\sLP
Around 1984 Macdonald visited his good friend Tonny Springer in Utrecht
and spoke in the staff colloquium there about his constant term conjectures
\cite{Macdonald 1982}.
The \emph{constant term conjecture} says that
\[ 
\mathrm{CT}(\Delta^{2k})=\prod_{j=1}^n\binom{kd_j}{k}\quad
\mbox{for all $k\in\mathbb{N}$.}
\]
Here $\Delta$ is the Weyl denominator,
$\mathrm{CT}(p)$ stands for the constant term
of the Fourier polynomial $p$, and $d_1,\ldots,d_n$ are
the degrees of the fundamental invariants of the Weyl group
for the root system~$R$ of rank~$n$ under consideration.
For the classical series
$\mathrm{A}_n,\mathrm{B}_n,\mathrm{C}_n,\mathrm{D}_n$
this formula was known by work of Selberg, Dyson and Good.
But for the exceptional root systems
$\mathrm{E}_6,\mathrm{E}_7,\mathrm{E}_8,\mathrm{F}_4,\mathrm{G}_2$
it still was an open question. 

During the years 1985--87 Gert Heckman and Eric Opdam
\cite{RSHF1}, \cite{RSHF2}, \cite{RSHF3}, \cite{RSHF4} developed a theory
of multivariable hypergeometric functions, depending on a (multi)parameter $k$,
which for special values of $k$ coincided with Harish-Chandra's spherical
functions. A crucial step in their work was to show the existence of a system
of commuting partial differential operators having as eigenfunctions
the functions they searched for.
So they were after a natural deformation in $k$
of the commutative algebra of PDOs found by Harish-Chandra in the group case.
For symmetric spaces of rank~1, including for instance hyperbolic spaces,
the spherical functions are indeed Euler--Gauss hypergeometric functions
(in the compact case Jacobi polynomials),
as Harish-Chandra (and \'Elie Cartan in the compact case) already showed.
There the extension to general $k$ is immediate.
For the root systems $\mathrm{BC}_2$ en $\mathrm{A}_2$
Tom Koornwinder \cite{Koornwinder 1974}, as part of his thesis in 1975,
had already given the commuting differential operators and corresponding
polynomial eigenfunctions. Also for general rank these matters had been treated
for root systems of type A (Jack polynomials) and BC, including explicit
expressions for the commuting PDOs (Sekiguchi, Debiard).
But it was the aim of Heckman and Opdam to give a uniform approach for all
root systems, not necessarily with explicit expressions.

In January 1988 Opdam defended in Leiden his thesis (essentially the papers
\cite{RSHF1}, \cite{RSHF3}, \cite{RSHF4}), with Macdonald present as
a committee member. There, in the final chapter \cite{RSHF4},
he succeeded to prove by analytic methods the existence of this 
hypergeometric system for general root systems.
This also implied the complete orthogonality of the
\emph{Jacobi polynomials associated with root systems}.
Shortly after Opdam showed
that, surprisingly, the results of his thesis could be applied to prove
the constant term conjecture for general root systems \cite{Opdam 1989}.

Almost parallel, in 1987, Macdonald had circulated handwritten manuscripts
which described the polynomials now known as Macdonald polynomials.
His first manuscript treated a class of symmetric polynomials with two
parameters $q$ and $t$. These extrapolated both the Jack polynomials and the
Hall--Littlewood polynomials by adding a parameter $q$.
He also gave an explicit system of commuting difference operators having
his polynomials as eigenfunctions. In the same year Simon
Ruijsenaars~\cite{Ruijsenaars 1987} had independently found this explicit
system. Macdonald gives an extensive treatment of this
class of symmetric polynomials in \cite[Ch.~VI]{Macdonald 1995}.
This chapter also includes
the so-called \emph{Pieri formula} \cite[\S VI.6]{Macdonald 1995}
for these polynomials, which Koornwinder \cite{Koornwinder 1988} proved in 1988
by showing a duality between the variables
and the spectral parameters. For this duality Koornwinder was inspired
by the two-variable case,
where the (homogeneous) Macdonald polynomials are essentially one-variable
continuous $q$-ultraspherical polynomials. There he could immediately observe
the duality, familiar as he was with the $q$-polynomial stuff by the work of
his former teacher Dick Askey \cite{Askey--Wilson 1985}.

Macdonald's next two manuscripts in 1987 concerned polynomials now known
as \emph{Macdonald polynomials associated with root systems}.
They are the direct $q$-analogues
of Jacobi polynomials associated with root systems as introduced by Heckman
and Opdam in the same time. Their polynomials could 
now be seen as $q\to1$ limits of Macdonald's polynomials.
For the root system of type $\textrm{A}$, Macdonald's $q$-polynomials are
contained in his class of symmetric functions,
so it was after the first manuscript a natural
step to pass to general root systems. Possibly the work of Heckman and Opdam,
which he had closely witnessed by several visits to Leiden in 1987
and 1988, has also been a motivation for him to look at $q$-analogues of their
polynomials.
The $q$-parameter in Macdonald's $q$-polynomials gave an additional
freedom by which he could prove full orthogonality in an elementary way.
And thus,
taking the limit for $q\to1$, the full orthogonality of the Heckman--Opdam
polynomials now also had an elementary proof.
Macdonald also gave conjectures for the norm and the evaluation at unity of
his polynomials.
In the limit case $q\to1$ his conjectures were in accordance with the
results of Heckman and Opdam.
In 1988 Macdonald put together the two manuscripts
into one manuscript \cite{Macdonald 1988}.
There he also mentioned Opdam's thesis and proof of the constant term
conjecture.

\begin{figure}[h]
\begin{subfigure}[t]{.4\linewidth}
\captionsetup{width=0.9\linewidth}
\includegraphics[height=4.8cm]{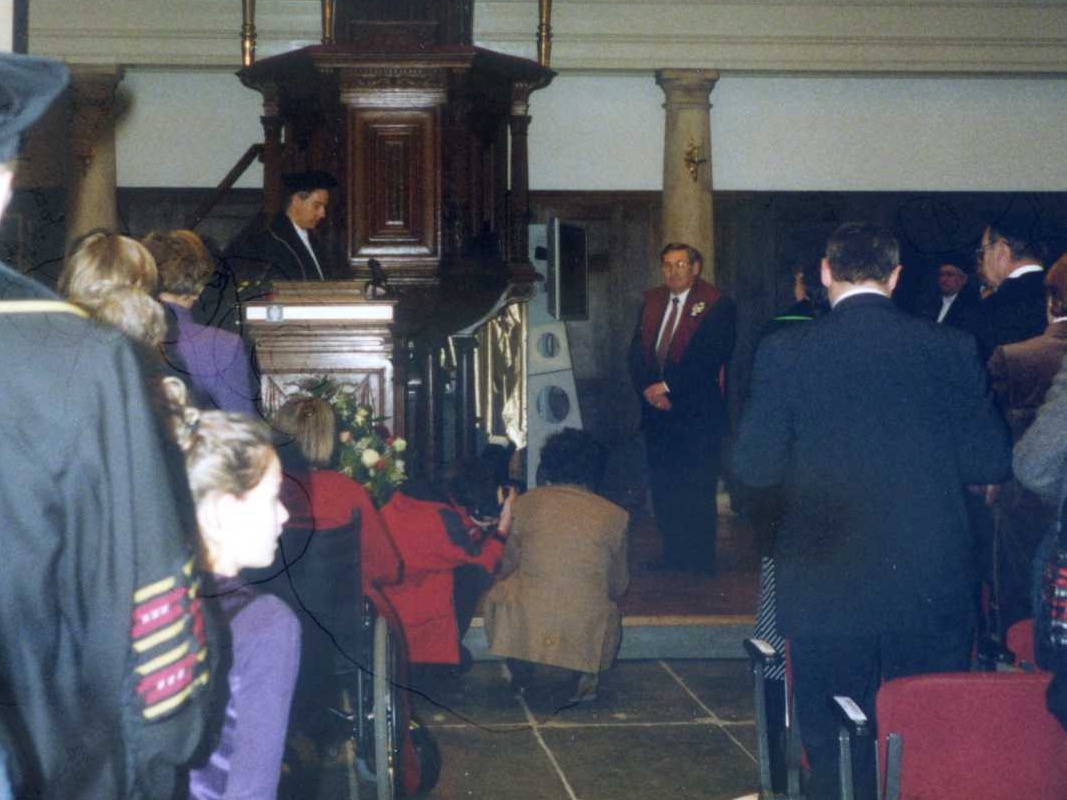}
\caption*{On January 8, 2002, in the aula of the University of Amsterdam,
Eric Opdam is reading the laudatio for the honorary doctorate of Macdonald
(standing to the right of the pulpit).}
\end{subfigure}\quad
\begin{subfigure}[t]{.5\linewidth}
\captionsetup{width=.9\linewidth}
\includegraphics[height=4.8cm]{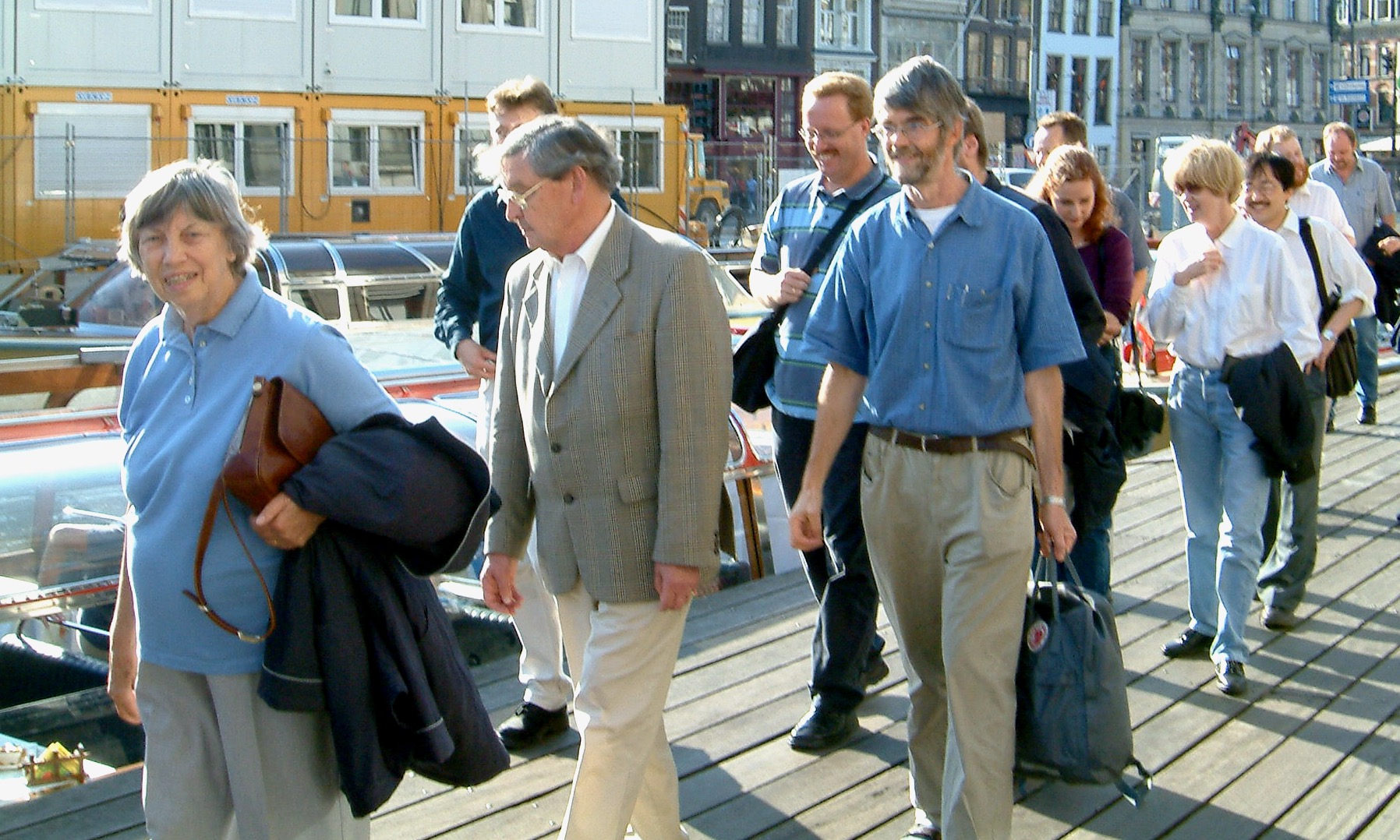}
\caption*{Greta and Ian Macdonald in front of attendees of a workshop on the
occasion of Tom Koornwinder's (middle) 60th birthday, all going to embark for a
roundtrip through the Amsterdam canals (August 2003).}
\end{subfigure}
\end{figure}

\paragraph{Further developments in connection with Macdonald polynomials}
\quad\sLP
In the spring of 1989 Heckman \cite{Heckman 1991} gave a
toric extension of Charles Dunkl's work
\cite{Dunkl 1989} on the differential-reflection operators which are now called
\emph{Dunkl operators}.
This enabled him to give an elementary construction of the hypergeometric
differential operators
which Opdam \cite{RSHF4} had constructed by transcendental methods,
and it led also to an elementary proof of the orthogonality
of their polynomial eigenfunctions, and of the norm
and evaluation at unity of the polynomials. A natural question was now how to
give a similar construction for Macdonald polynomials.
This was quickly answered.

\begin{wrapfigure}{R}{0.35\textwidth}
\centering
\captionsetup{width=0.9\linewidth}
\includegraphics[width=0.35\textwidth]{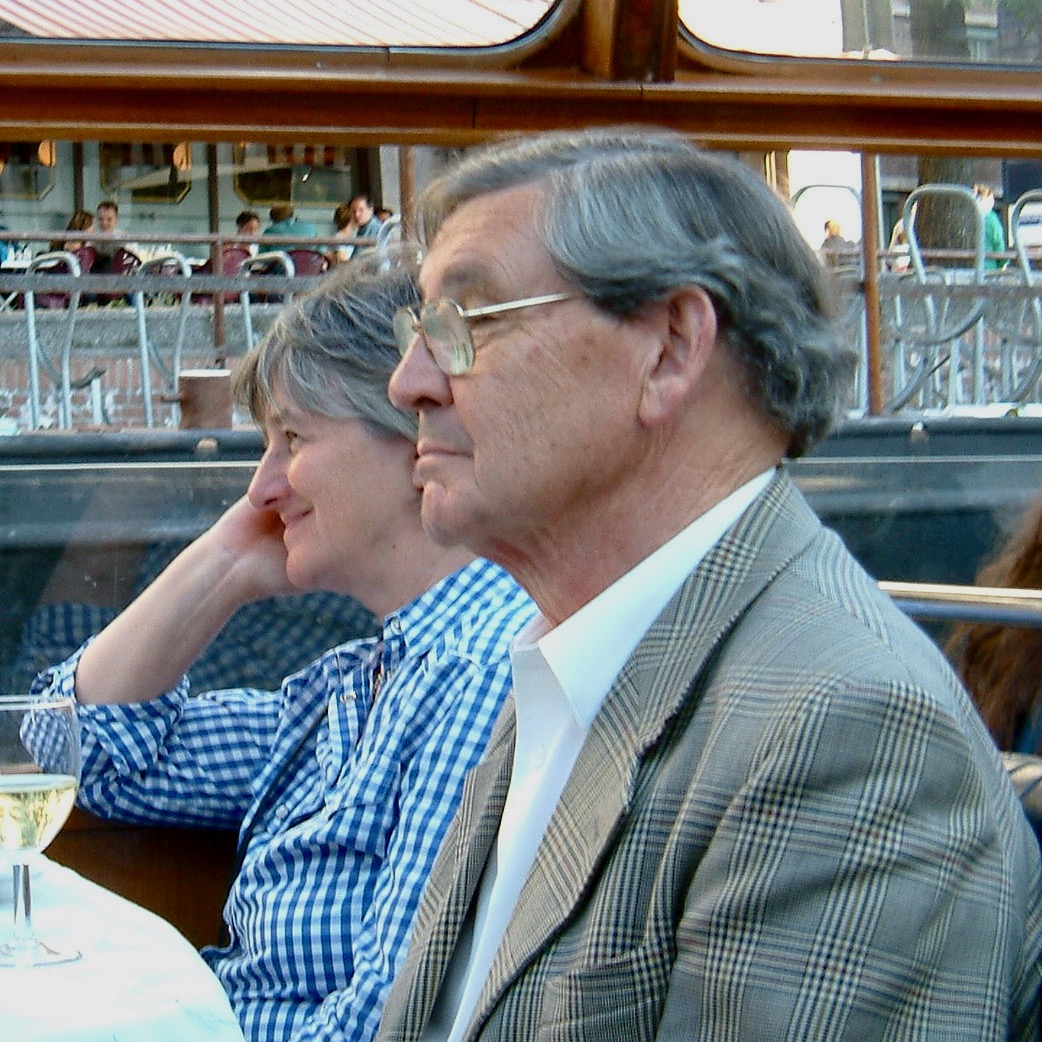}
\caption*{Macdonald during the roundtrip through the canals at the mentioned
workshop in Amsterdam, August 2003. Beside him Suzanne Faraut.}
\end{wrapfigure}

In the fall of 1989 the Berlin wall was torn down, and in the summer of the 
next year many Russians visited the ICM in Kyoto, among whom
there was Ivan Cherednik. 
In a few years all conjectures of Macdonald were made theorems by Cherednik 
\cite{Cherednik 1995}. 
It was not always easy reading in Ivan's papers, but expositions such as by
Opdam \cite{Opdam 1995} and Macdonald \cite{Macdonald 2003} were helpful for
understanding. The key to all this were \emph{Hecke algebras}:
first the affine Hecke algebra of Iwahori and Matsumoto coming from the
$p$-adic world, next the degenerate 
Hecke algebra of Drinfeld \cite{Drinfeld 1986}, and for the duality
between variable and spectral parameter the 
\emph{double affine Hecke algebra} of Cherednik.
The commuting generalized Dunkl operators had \emph{non-symmetric}
(non-$W$-invariant) Macdonald polynomials as eigenfunctions, and these were
projected onto the symmetric ($W$-invariant) Macdonald polynomials by
the symmetrization map.

Another new development were the \emph{Koornwinder polynomials}
\cite{Koornwinder 1992}. These are a significant 5-parameter extension
of the 3-parameter family of Macdonald polynomials
associated with the non-reduced root system of type
$\mathrm{BC}_n$ (here $q$ is not considered as a parameter).
This was at the same time an $n$-variable extension
of the $4$-parameter family of Askey--Wilson-polynomials
\cite{Askey--Wilson 1985} in one variable.
Also for this case Macdonald formulated conjectures, which again could be
proved by Hecke algebra methods \cite{Sahi 1999}, \cite{Macdonald 2003}.
Younger Dutch mathematicians were stimulated by these developments, see
van Diejen \cite{Diejen 1996} and Stokman \cite{Stokman 2000} (former students
of Ruijsenaars and Koornwinder, respectively). Earlier van der Lek
wrote an influential thesis \cite{Lek 1983}. This was helpful for Ion
\cite{Ion 2003} to give a conceptual proof of the duality theorem for
double affine Hecke algebras which was proved in Macdonald's last book
\cite[Theorem (3.5.1)]{Macdonald 2003} using a case-by-case analysis.

In this period, starting about 1985, the \emph{quantum groups} also entered
the stage. Gradually it became clear that $q=\mbox{quantum}$:
Macdonald polynomials and (to a lesser extent)
Koornwinder polynomials
lived for special parameter values as spherical functions on quantum groups
\cite{Letzter 2004},~\cite{Noumi 1997}.

The authors of this paper, each from his own
perspective, have experienced the described period as extremely fascinating.
Many developments, motivated in quite different ways, interacted and
strengthened each other.
\paragraph{Macdonald as a person}\quad\sLP
We conclude with some words about Ian Macdonald as a person and as a
mathematician. He was a very gifted speaker.
During his lecture in S\'eminaire Bourbaki in the mid 1990s
\cite{Macdonald 1996} 
one of us was sitting beside Jim Arthur, an expert in the
\emph{Langlands program}. At the end of the lecture he was silent for a while,
and then said: ``What an outstanding expositor''. 
We felt indeed privileged that Macdonald visited the Netherlands so frequently
to share his views with us.
On these occasions Ian was usually accompanied by his wife Greta, originating
from the Flemish part of Belgium. Ian had no driving license, and thus she
drove while he read the map for navigation. During a lecture in the early 1990s
Ian started as follows:
\newpage\indent
If you have a mathematical problem, there are four things you can do: 
\sLP\indent\indent
1. Ignore it. If you can't: \\\indent\indent
2. Prove it. If you can't:\\\indent\indent
3. Disprove it. Finally if you can't do that either:\\\indent\indent
4. Generalise it.
\bLP
Everybody was laughing, but of course he was also very serious.
If you have already considered options 2 and 3 in vain, then it would be a pity
to throw the problem away before having tried option 4. 
He seemed to apologize that in his lecture he would look again more broadly
than asked for by the problem. His conjectured generalizations always turned
out to be correct, and in this way he greatly stimulated others.

No doubt Ian was a mathematician with a very classical inclination.
During his education in Cambridge, the spirit of Hardy and Ramanujan may
still have been felt. Ian was a master in dealing with formulas, but in
retrospect all these formulas were directed to a better understanding
of the  representation theory of the simple algebraic groups
(over a field or ring).
We are truly grateful that the mathematical road travelled by Ian in his later
life crossed so frequently with ours.

\paragraph{Acknowledgements}
The first two photos are by Peter van Emde Boas, and the last two ones
by Margrit R\"osler and Michael Voit. We thank the referee for suggesting
various textual improvements.

\noindent
Gert Heckman, Radboud University Nijmegen, \url{g.heckman@math.ru.nl} \\
Tom Koornwinder, University of Amsterdam, \url{thkmath@xs4all.nl} \\
Eric Opdam, University of Amsterdam, \url{e.m.opdam@uva.nl}

\end{document}